# ON THE INTEGRAL OF GEOMETRIC BROWNIAN MOTION

Michael Schröder

**Abstract.** This paper studies the law of any power of the integral of geometric Brownian motion over any finite time interval. As its main results, two integral representations for this law are derived. This is by enhancing the Laplace transform ansatz of [**Y**] with complex analytic methods, which is the main methodological contribution of the paper. The one of our integrals has a similar structure to that obtained in [**Y**], while the other is in terms of Hermite functions as those of [**Du01**]. Performing or not performing a certain Girsanov transformation is identified as the source of these two forms of the laws. For exponents equal to 1 our results specialize to those obtained in [**Y**], but for exponents equal to minus 1 they give representations for the laws which are markedly different from those obtained in [**Du01**].

Keywords. Brownian motion; Laplace transform; Asian options.

AMS 2000 subject classification. 60J65; 91B28.

**1. Introduction:** This paper works towards a synthesis of the representations derived in the recent years for the law of the integral of geometric Brownian motion over a finite time interval. These processes have both a surprisingly rich theory and manifold applications ranging from the physics of random media to mathematical finance and insurance. In fact, the insurance motivated study of certain perpetuities in [**Du90**] seems to have initiated this line of research. Here the integral of geometric Brownian motion over the whole time axis is considered and is shown to be distributed as the reciprocal of a certain gamma variable In mathematical finance, on the other hand, valuing the so called Asian options asks for an as explicit as possible expression for a certain functional of the integral of geometric Brownian motion over a finite time interval. Drawing on his probabilistic interpretation of the Hartman–Watson identities in [**Y80**], Yor was able to determine the law of these processes in [**Y**]. His approach, based on the Laplace transform and using Bessel processes, has proved both fertile and fundamental, and the string of papers he wrote about it are now collected in [**Y01**]. In a recent developement initiated in [**Du00**], however, the focus has shifted towards studying the reciprocal of the integral of geometric Brownian motion. As one of the lastest results, [**Du01**] derives an integral for its law in terms of Hermite functions. This is a structural advance in that Hermite functions form a well–studied class of special functions which generalize both the error function and the Hermite polynomials and which come from boundary value problems in potential theory for domains whose surface is an infinite parabolic cylinder. They have already appeared in [**SA**] in connection with valuing Asian options, and this suggests a Hermite function form of the law of the integral of geometric Brownian motion itself.

This paper now has the more general aim so to develop a common framework for Yor's and Dufresne's results by focussing on the structure of the laws of the integral of geometric Brownian motion raised to any real power $\varepsilon$. This seems to be different in spirit from the work which aims at a better understanding of this circle of ideas by developing alternatives to Yor's approach like those of [**AMS**] based on the Feynman–Kac ansatz. In our framework, specializing $\varepsilon$ to 1 gives results about the law of the integral of geometric Brownian motion, and specializing $\varepsilon$ to minus 1 gives results about the laws of the reciprocals of these processes. As the main findings of this paper, two forms of these $\varepsilon$–power laws are discussed in Part I: an integral representation of a type similar to the density obtained in [**Y**] and a second in terms of Hermite functions. The first of these specialize in particular to the representations for the law of of the integral of geometric Brownian motion resulting from [**Y**]. Specializing to the laws of the reciprocals of these processes, however, the resulting Hermite function forms are markedly different from those obtained in [**Du01**], and the exact relations between these two representations remain to be studied. And here we should add that mutatis mutandis this also holds for the Barnes type integral representations in terms of Whittaker functions of the laws of the integral of geometric Brownian motion which figure in [**CMY**, §4.1].

Proved are our Yor type integrals for the $\varepsilon$–power laws in Part III by a generalization of Yor's original approach based on the Laplace transform and using Bessel processes. However, it is the approach developed for establishing the Hermite function form of these laws which seems to be the main methodological contribution of this paper. As explained in Parts IV to VI this is by enhancing the generalized Yor approach of Part III with complex analytic methods. In fact, the Yor form of the laws is found to result on performing a Girsanov transformation which reduces the drift terms of the geometric Brownian motion to zero. Not performing such a Girsanov transformation one gets the Hermite function form of the laws. In fact, they result directly if the drift coefficients are positive and by using an analytic continuation procedure otherwise. So the following philosophy seems to be tempting: an application of a Girsanov transformation can be replaced by an introduction of certain complex analytic methods and vice versa. The discussion of Part VI, however, shows that the situation is not that straightforward. These two methodologies are not equivalent but complement each other in a rather intricate way. Still, we consider this as another example why combining stochastics with complex analytic methods seems to be a rather promising line of strategy and study.

## Part I  Statement of results

**2. Basic concepts:** The original motivation for this aper was a description of the densities of Yor's accumulation processes $A^{(\nu)}$, the *Asia densities* $\alpha^{(\nu)}$, and the densities of the reciprocals of these processes, the *reciprocal Asia densities* $\beta^{(\nu)}$. Results about these two types of densities, however, are implied by those about the densities of any $\varepsilon$–th power of Yor's accumulation processes $A^{(\nu)}$, for $\varepsilon$ any non–zero real. They are called $\varepsilon$-*power Asia densities*, are denoted by $\alpha^{(\nu,\varepsilon)}$, and form the actual object of study of this chapter. By definition, to any random variable

$$\left(A_t^{(\nu)}\right)^\varepsilon \qquad \text{where} \qquad A_t^{(\nu)} = \int_0^t e^{2(\nu x + B_x)}\, dx$$

is thus associated its density $w \mapsto \alpha_t^{(\nu,\varepsilon)}(w)$ as a map on the non–negative real line, and we have by construction $\alpha^{(\nu)} = \alpha^{(\nu,1)}$ and $\beta^{(\nu)} = \alpha^{(\nu,-1)}$. We are looking for integral representations of these densities. Their principal structure is that as the product of the functions $c_{\nu,\varepsilon,t}$ given by

$$c_{\nu,\varepsilon,t}(w) = |\varepsilon|^{-1} \frac{2^{\frac{\nu}{2}}}{\sqrt{\pi t}} \exp\left(-\frac{1}{2}\left(\nu^2 t + w^{-\frac{1}{\varepsilon}}\right)\right) w^{\frac{1}{2\varepsilon}(\nu-1)-1}$$

times integrals over the contours of integration $\log C_{\theta,R}$ given by

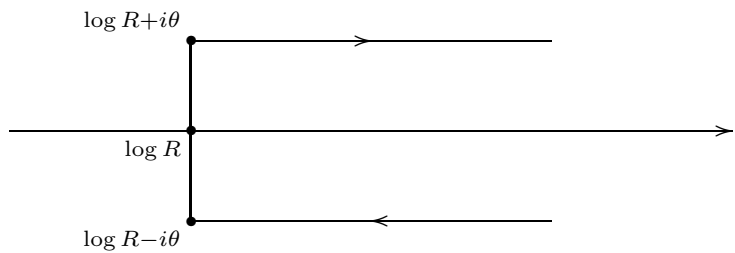

**Figure 1.** The contour $\log C_{\theta,R}$.

as the logarithms of the Hankel contours $C_{\theta,R}$ with angles $\theta$ in $[\frac{\pi}{2}, \pi]$ and radii $R \geq 1$:

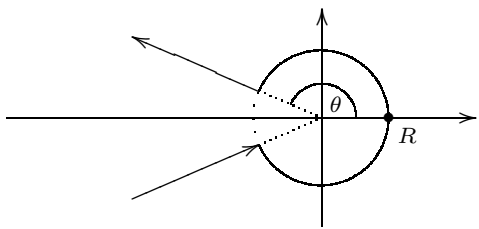

**Figure 2.** The Hankel contour $C_{\theta,R}$.

We study two forms of the $\varepsilon$–power Asia densities. One in the spirit of the results of [**Y**] for Asia densities and to other in the spirit of [**SH**] or [**Du01**] using Hermite functions. Following [**L**, §§10.2ff], these functions are given in terms of the Kummer confluent hypergeometric function $\Phi$. The *Hermite function* $H_\mu$ of degree any complex number $\mu$ is the function on the complex plane defined by

$$H_\mu(z) = \frac{2^\mu \cdot \Gamma(1/2)}{\Gamma((1-\mu)/2)} \cdot \Phi\left(-\frac{\mu}{2}, \frac{1}{2}, z^2\right) + z \cdot \frac{2^\mu \cdot \Gamma(-1/2)}{\Gamma(-\mu/2)} \cdot \Phi\left(\frac{1-\mu}{2}, \frac{3}{2}, z^2\right).$$

Hermite functions $H_\mu$ so are holomorphic on the complex plane as functions of both their variable $z$ and their degree $\mu$. If the real part $\operatorname{Re}(\mu)$ of $\mu$ is bigger than $-1$, they have the integral representation

$$H_\mu(z) = \frac{2^{\mu+1}}{\sqrt{\pi}} e^{z^2} \int_0^\infty e^{-x^2} x^\mu \cos\left(2zx - \frac{1}{2}\mu\pi\right) dx$$

and so specialize to the $\mu$–th Hermite polynomials if $\mu$ is any non–negative integer. If the real part of $\mu$ is negative, they have the integral representation:

$$H_\mu(z) = \frac{1}{\Gamma(-\mu)} \int_0^\infty e^{-u^2 - 2zu} u^{-(\mu+1)} du,$$

and so specialize via $(2/\sqrt{\pi})H_{-1}(z) = \exp(z^2)\operatorname{Erfc}(z)$ to the complementary error function Erfc. They are moreover connected with the parabolic cylinder functions $D_\mu$ and with the Kummer confluent hypergeometric function of the second kind $\Psi$.

**3. Yor type $\varepsilon$–power Asia densities:** This section discusses generalizations of Yor's original form of the Asia densities resulting from [**Y**] to $\varepsilon$–power Asia densities. In terms of the concepts of §2, the basic result is the following contour integral representation.

**Theorem:** *For any real $\nu$ and non–zero real $\varepsilon$, the $\varepsilon$–power Asia density $\alpha^{(\nu,\varepsilon)}$ is given by the following double integral*

$$\alpha_t^{(\nu,\varepsilon)}(w) = c_{\nu,\varepsilon,t}(w) \int_0^\infty x^\nu e^{-x^2} \frac{1}{2\pi i} \int_{\log C_{\theta,R}} e^{-\frac{\xi^2}{2t}} \sinh(\xi) \exp\left(\frac{2x\cosh(\xi)}{\sqrt{2w^{1/\varepsilon}}}\right) d\xi \, dx$$

*for any positive reals $t$ and $w$, and any real parameters $\theta$ in $[\frac{\pi}{2}, \pi]$ and $R \geq 1$.*

We discuss two specializations of the Theorem. If $\theta = \frac{\pi}{2}$ and $R = 1$, we obtain

$$\alpha_t^{(\nu,\varepsilon)}(w) = \frac{c_{\nu,\varepsilon,t}(w)}{\pi} e^{\frac{\pi^2}{8t}} \int_0^\infty x^\nu e^{-x^2} \int_0^\infty e^{-\frac{y^2}{2t}} \cosh(y) \cos\left(\frac{2x}{\sqrt{2w^{1/\varepsilon}}} \sinh(y) - \frac{\pi y}{2t}\right) dy \, dx$$

for any positive reals $t$ and $w$. In the case $\theta = \pi$ and $R = 1$ as in [**Y**], we get

$$\alpha_t^{(\nu,\varepsilon)}(w) = \frac{c_{\nu,\varepsilon,t}(w)}{\pi} e^{\frac{\pi^2}{2t}} \int_0^\infty x^\nu e^{-x^2} \int_0^\infty e^{-\frac{y^2}{2t}} \sinh(y) \sin\left(\frac{\pi y}{t}\right) \exp\left(-\frac{2x\cosh(y)}{\sqrt{2w^{1/\varepsilon}}}\right) dy \, dx$$

for any positive reals $t$ and $w$. For $\varepsilon = 1$ this last formula is seen to specialize to the double integral for the Asia density resulting from [**Y**, §6].

**4. The Hermite function form of $\varepsilon$–power Asia densities:** In the circle of ideas which has developed around understanding accumulation processes $A^{(\nu)}$, Hermite functions as recalled in §2 seem to have been explicitly identified first in connection with Asian options. Here they appear in the integral representations of [**SA**]. Subsequently, they were shown to play a role for reciprocal Asia densities $\beta^{(\nu)}$ in [**Du01**]. This could be surprising given in particular Yor's results in [**Y**] on Asia densities as we have generalized them in §3. However, it should be not. In fact, this paper offers what seems to be a rather natural explanation for the appearance of these functions. Using the concepts of §2, the precise result is the following representation of $\varepsilon$–power Asia densities as contour integrals in terms of Hermite functions $H_\mu$.

**Theorem:** *For any real $\nu$ and any non–zero real $\varepsilon$, the $\varepsilon$–power Asia density $\alpha^{(\nu,\varepsilon)}$ is given by the following contour integral*

$$\alpha_t^{(\nu,\varepsilon)}(w) = \Gamma(\nu+1) \frac{c_{\nu,\varepsilon,t}(w)}{2\pi i} \int_{\log C_{\theta,R}} e^{-\frac{\xi^2}{2t}} \sinh(\xi) \, H_{-(\nu+1)}\left(-\frac{\cosh(\xi)}{\sqrt{2w^{1/\varepsilon}}}\right) d\xi$$

*for any positive reals $t$ and $w$.*

**Remark:** If $\nu$ is any negative integer, the density is given by the product of $c_{\nu,\varepsilon,t}$ times the product of the residue $(-1)^n/n!$ of the gamma function in $n = -(\nu+1)$ and the corresponding contour integral in terms of the derivative of $H_{-\mu}$ with respect to $\mu$ at $\mu = \nu+1$. This last function is a new higher transcendental function which deserves further study.

We discuss two specializations of the Theorem. If $\theta = \frac{\pi}{2}$ and $R = 1$, we obtain

$$\alpha_t^{(\nu,\varepsilon)}(w) = \Gamma(\nu+1)\frac{c_{\nu,\varepsilon,t}(w)}{\pi} e^{\frac{\pi^2}{8t}} \int_0^\infty e^{-\frac{y^2}{2t}} \cosh(y) \operatorname{Re}\left(e^{-\frac{\pi y i}{2t}} H_{-(\nu+1)}\left(-i\frac{\sinh(y)}{\sqrt{2w^{1/\varepsilon}}}\right)\right) dy$$

for any positive reals $t$ and $w$. In the case $\theta = \pi$ and $R = 1$, we get

$$\alpha_t^{(\nu,\varepsilon)}(w) = \Gamma(\nu+1)\frac{c_{\nu,\varepsilon,t}(w)}{\pi} e^{\frac{\pi^2}{2t}} \int_0^\infty e^{-\frac{y^2}{2t}} \sinh(y) \sin\left(\frac{\pi y}{t}\right) H_{-(\nu+1)}\left(\frac{\cosh(y)}{\sqrt{2w^{1/\varepsilon}}}\right) dy$$

for any positive reals $t$ and $w$. If $\nu > -1$, these results visibly coincide with those of §3. For $\varepsilon = -1$ they give representations for the reciprocal Asia densities $\beta^{(\nu)}$ which are similar to but markedly different from those of [**Du01**].

## Part II    Three basic auxiliary results

**5.  A basic Laplace transform:**   A key result of the Laplace transform approach to analyzing problems involving Yor's accumulation processes

$$A_t^{(\mu)} = \int_0^t e^{2(B_w + \mu w)} dw \qquad t \in [0, \infty)$$

has been the factorization of geometric Brownian motion afforded by the Lamperti identity

$$e^{B_t + \mu t} = R^{(\mu)}\left(A_t^{(\mu)}\right),$$

proved for instance in [**SH**, Chapter II §5]. Herein $R^{(\mu)}$ is the Bessel process of index $\nu$ started at 1 at time $t = 0$. We have distilled from [**Y**] the following consequence of this factorization for computing Laplace transforms of expectations in Yor's accumulation process which is a basic tool for computing the above forms of the $\varepsilon$–power Asia densities.

**Proposition:**   *Let $\varepsilon$ be any non–zero real. If $\mu \geq 0$, we have for any measurable functions $f$, $g$ on the non–negative real line*

$$\int_0^\infty e^{-zt} E\left[f\left((A_t^{(\mu)})^\varepsilon\right) g\left(e^{B_t}\right)\right] dt$$

$$= \int_0^\infty \int_0^\infty f(y) g(\rho) |\varepsilon|^{-1} \frac{\rho^{\mu-1}}{y} \exp\left(-\frac{1+\rho^2}{2y^{1/\varepsilon}}\right) I_{\sqrt{2z+\mu^2}}\left(\frac{\rho}{y^{1/\varepsilon}}\right) d\rho \, dy$$

*for any real $z$ with a sufficiently big positive real part.*

This is based on the following time change result of independent interest

**Lemma:** *If $\mu \geq 0$, we have for any measurable function $F$ on the non–negative real line*

$$E\left[\int_0^\infty e^{-zt} F(A_t^{(\mu)}) \, dt\right] = E\left[\int_0^\infty e^{-z\tau^{(\mu)}(w)} F(w) \frac{dw}{(R^{(\mu)})^2(w)}\right]$$

where $\tau^{(\mu)}(w) = \inf\{u \, | \, A^{(\mu)}(u) > w\}$.

**Proof of the Lemma:** For the proof notice that by construction $w \mapsto \tau^{(\mu)}(w)$ is the inverse map of $t \mapsto A^{(\mu)}(t)$. Time changing thus gives

$$E\left[\int_0^\infty e^{-zt} F(A_t^{(\mu)}) \, dt\right] = E\left[\int_0^\infty e^{-z\tau(w)} F(w) \, d\tau^{(\mu)}(w)\right].$$

In fact, $A^{(\mu)}$ starts at zero at time $t = 0$ by construction, and the strong law of large numbers for Brownian motion implies that $A^{(\mu)}(t)$ goes to infinity with $t$ to infinity. To interpret the differential, apply the inverse function theorem to get

$$d\tau^{(\mu)}(w) = \frac{dw}{\exp\left(2(\mu x + B_x)\right)\big|_{x=\tau^{(\mu)}(w)}}.$$

Using the above Lamperti identity, the denominator is equal to the square of the Bessel process of index $\mu$ at time $w$, as was to be shown for the Lemma.

**Proof of the Proposition:** There are two main ideas for establishing the Proposition by explicitly computing the Laplace transform

$$L(z) = \int_0^\infty e^{-zt} E\left[f\left((A_t^{(\mu)})^\varepsilon\right) g\left(e^{B_t}\right)\right] dt.$$

The first is to make the simpler process $\exp(B_w)$ more complicated in a controlled fashion using Lamperti's identity. Moreover using Tonelli's theorem this gives

$$L(z) = E\left[\int_0^\infty e^{-zt} f\left((A_t^{(\mu)})^\varepsilon\right) g\left(R^{(\mu)}(A_t^{(\mu)})\right) dt\right].$$

Applying the time change result of the above Lemma, also based on Lamperti's identity,

$$L(z) = E\left[\int_0^\infty e^{-z\tau^{(\mu)}(w)} f(w^\varepsilon) \, g(R_w^{(\mu)}) \frac{dw}{(R_w^{(\mu)})^2}\right].$$

Inside the expectation change variables $w = y^{1/\varepsilon}$. Then reverse the Tonelli theorem and integrate against the semigroup densities of the resulting index-$\mu$ Bessel processes to get

$$L(z) = \int_0^\infty \int_0^\infty E\left[e^{-z\tau^{(\mu)}(y^{1/\varepsilon})} \Big| R_{y^{1/\varepsilon}}^{(\mu)} = \rho\right] \frac{f(y)g(\rho)}{\rho^2} p_{y^{1/\varepsilon}}^{(\mu)}(1,\rho) \, d\rho \, \frac{y^{\frac{1}{\varepsilon}-1}}{|\varepsilon|} \, dy.$$

Here $p_w^{(\mu)}$ is the time–$w$ semigroup of the Bessel process of index $\mu$ starting at 1 at time $t = 0$. Since $\mu$ is non–negative, from [**RY**, X §1] recall

$$p_w^{(\mu)}(1,\rho) = \frac{\rho^{\mu+1}}{w} \exp\left(-\frac{1+\rho^2}{2w}\right) I_\mu\left(\frac{\rho}{w}\right).$$

The second idea is to handle the conditional expectation factor in the integrand of the last double integral using Yor's Bessel expectation [**Y80**, Théorème 4.7, p.80]. Recall this result asserts, for any non–negative real $\mu$, that we have

$$E\left[e^{-z\cdot \tau_w^{(\mu)}} \Big| R_w^{(\mu)} = x\right] = \frac{I_{\sqrt{2z+\mu^2}}}{I_\mu}\left(\frac{x}{w}\right),$$

for any positive reals $z$, $w$ and $x$. On substitution the proof of the Proposition is complete.

**6. A basic Laplace inverse:** This section explicitly reviews how the Bessel function factors of the integrands of the identity of §5 Proposition are obtained as Laplace transforms. This is already at the base of the inversion argument of [**Y80**, §5], and using the concepts of §2, the precise result is the following

**Lemma:** *For any complex numbers $\eta$ with positive real part and $\mu$,*

$$\mathscr{L}^{-1}\left(I_{\sqrt{2z+\mu^2}}(\eta)\right)(t) = \eta \frac{e^{-\frac{1}{2}\mu^2 t}}{\sqrt{2\pi t}} \frac{1}{2\pi i} \int_{\log C_{\theta,R}} e^{-\frac{\xi^2}{2t} + \eta \cosh(\xi)} \sinh(\xi)\, d\xi$$

*for any reals $t > 0$, $\theta$ in $[\frac{\pi}{2}, \pi]$ and $R \geq 1$.*

The proof of the Lemma is based on the following classical Hankel contour integral for the modified Bessel function $I$ which can be found in [**WW**, 17·231, p.362], for instance:

$$I_\rho(\eta) = \frac{1}{2\pi i}\int_{\log C_{\theta,R}} e^{-\rho\cdot\xi + \eta\cdot\cosh(\xi)}\, d\xi,$$

for any complex numbers $\eta$ with positive real part and $\rho$. To identify the first exponential factor of this last integral's integrand as a Laplace transform if $\rho = (2z+\mu^2)^{1/2}$, recall from [**D**, Beispiel 8, p.50f] the standard Laplace transform

$$e^{-a\sqrt{z}} = \mathscr{L}(\psi_a)(z) \qquad \text{where} \qquad \psi_a(t) = \frac{a}{2\sqrt{\pi}\, t^3} e^{-\frac{a^2}{4t}}$$

for any complex numbers $z$ with positive real part and $a$ such that both $a$ and its square $a^2$ have positive real parts. If $a$ is in the contour $\log C_{\theta,R}$, this last condition can be checked to be satisfied if the radius $R$ is large enough. First considering this case, we then have

$$e^{-\xi\sqrt{2z+\mu^2}} = \int_0^\infty e^{-zt} e^{-\frac{1}{2}\mu^2 t} \frac{\eta}{\sqrt{2\pi t^3}} e^{-\frac{\xi^2}{2t}}\, dt,$$

for any complex $z$ with sufficiently big positive real part, whence on substitution

$$I_{\sqrt{2z}}(\eta) = \frac{e^{-\frac{1}{2}\mu^2 t}}{\sqrt{2\pi t^3}} \frac{1}{2\pi i} \int_{\log C_{\theta,R}} e^{\eta \cdot \cosh(\xi)} \int_0^\infty e^{-zt} \frac{e^{-\frac{1}{2}\mu^2 t}}{\sqrt{2\pi t^3}} \xi e^{-\frac{\xi^2}{2t}} \, dt \, d\xi.$$

The point now is to justify in this double integral interchanging the order of the Laplace transform with the integration over the Hankel contour. Moreover choosing $\theta$ in $(\pi/2, \pi]$ for this, the hyperbolic cosine in the second exponent has a negative real part. If $\eta > 0$ in addition, the integrand of the above double integral can be checked to go to zero in absolute value with $(z, \xi)$ going to infinity. Thus, under these restrictions, Fubini's theorem gives

$$\mathscr{L}^{-1}\left(I_{\sqrt{2z+\mu^2}}(\eta)\right)(t) = \frac{e^{-\frac{1}{2}\mu^2 t}}{\sqrt{2\pi t^3}} G_\eta(t)$$

setting

$$G_\eta(t) = \frac{1}{2\pi i} \int_{\log C_{\theta,R}} e^{\eta \cdot \cosh(\xi)} \xi e^{-\frac{\xi^2}{2t}} \, d\xi.$$

With its integrand entire, using Cauchy's theorem $G_\eta$ is independent of any chosen contour $\log C_{\theta,R}$ with $\theta \in [\pi/2, \pi]$ and $R \geq 1$, and partial integration then is seen to give

$$G_\eta(t) = \frac{\eta t}{2\pi i} \int_{\log C_{\theta,R}} e^{-\frac{\xi^2}{2t} + \eta \cosh(\xi)} \sinh(\xi) \, d\xi.$$

This completes the proof of the Lemma.

**7. A basic analyticity criterion:** This section addresses the construction of analytic functions by integration from complex valued functions $f$ on product spaces of the form $V \times X$. Here $V$ is any subset of the complex plane with a non–empty interior and $X$ is any measurable subset of a measure space with measure $\mu_X$. The precise result to be proved in this set up is the following *analyticity criterion*

**Lemma:** *For any integrable complex valued map $f$ on $V \times X$, the map on $V$ defined by*

$$F(v) = \int_X f(v, x) \, d\mu_X(x)$$

*is analytic on $V$ if for almost all $x$ in $X$ the map $v \mapsto f(v, x)$ is analytic on $V$.*

**Remark:** Integrability of any complex valued function $f$ on $V \times X$ is assured if there is an integrable map $g$ on $X$ such that $|f(v, x)| \leq g(x)$ for any $x$ in $X$.

With the Remark being clear, the idea for the proof of the Proposition is to establish analyticity of $F$ on $V$ by applying Morera's theorem, see [**R**, 10.17, p.208] for instance. For this we have to show

$$\int_\Delta F(v) \, dv = 0$$

for any triangle $\Delta$ contained in the interior of $V$. Writing out $F$ in terms of its defining integral, notice that the integrand $f$ of the resulting double integral is integrable by assumption. An application of Fubini's theorem thus gives

$$\int_\Delta F(v)\,dv = \int_0^\infty \int_\Delta f(v,x)\,dv\,d\mu_X(x)\,.$$

By assumption the integrands of almost all inner integrals are analytic on $V$. To these Cauchy's theorem as in [**R**, 10.13, p.205] applies and shows that they are equal to zero. Thus the integral of $F$ over $\Delta$ is equal to zero, as was to be shown.

# Part III
# Proof of the Yor type $\varepsilon$–power Asia densities

**8. Reductions for the Yor type $\varepsilon$–power Asia densities:** This section reviews following [**Y**] how computing the $\varepsilon$–power Asia density $\alpha^{(\nu,\varepsilon)}$, for any non–zero real $\varepsilon$, is reduced to computing conditional $\varepsilon$–power Asia densities $a_{\varepsilon,t}(x,w)$. This is based on the Girsanov transformation of measure such that $W_t = \nu t + B_t$ becomes a standard Brownian motion. Dropping reference to this new measure, we then have for any measurable function $f$ on the non–negative real line into intself

$$E^Q\left[f\left((A_t^{(\nu)})^\varepsilon\right)\right] = E\left[f\left((A_t^{(0)})^\varepsilon\right)g_\nu\left(e^{W_t}\right)\right] \cdot e^{-\frac{t}{2}\nu^2}$$

with maps $g_\nu(x) = x^\nu$. Computing this expectation so becomes a special case of computing expectations of the form

$$E\left[f\left((A_t^{(0)})^\varepsilon\right)g\left(e^{B_t}\right)\right]$$

with Borel maps $f$, $g$ from the non–negative reals into themselves. Getting the $\varepsilon$–*power Asia densities*, the densities of the $\varepsilon$–th power of $A^{(\nu)}(t)$, so is reduced to getting the *conditional $\varepsilon$–power Asia densities* which for any reals $t > 0$ and $x$ are given by the maps

$$w \mapsto a_{\varepsilon,t}(x,w) := Q\left((A_t^{(0)})^\varepsilon \in dw \big| B_t = x\right)$$

and which as densities of powers of $A^{(0)}(t)$ conditional on $B_t = x$ are independent of $\nu$.

**9. Laplace transforms of Yor type $\varepsilon$–power Asia densities:** This section adapts the computations of [**Y**] to computing the Laplace transform with respect to time of the conditional $\varepsilon$–power Asia densities

$$w \mapsto a_{\varepsilon,t}(x,w) := Q\left((A_t^{(0)})^\varepsilon \in dw \big| B_t = x\right)$$

defined in §8. The precise result in terms of modified Bessel functions $I_\mu$ is the following

**Lemma:** *For any non–zero real $\varepsilon$ and any reals $w > 0$, $x$, we have*

$$\frac{1}{\sqrt{2\pi t}} \exp\left(-\frac{x^2}{2t}\right) a_{\varepsilon,t}(x,w) = \frac{|\varepsilon|^{-1}}{w} \exp\left(-\frac{1+e^{2x}}{2w^{1/\varepsilon}}\right) \mathscr{L}^{-1}\left(I_{\sqrt{2z}}\left(\frac{e^x}{w^{1/\varepsilon}}\right)\right)(t).$$

*for any positive real $t$.*

The idea of the proof is to compute in two ways the Laplace transform

$$L(z) = \int_0^\infty e^{-zt} E\left[g(e^{B_t}) f(A_t^{(0)})\right] dt$$

for any Borel functions $f$, $g$ on the non–negative real line into itself and with $z \geq z_0$ any sufficiently big positive real. First, §5 Proposition applied with $\mu = 0$ gives on changing variables $\rho = \exp(x)$,

$$L(z) = \int_0^\infty \int_{\mathbf{R}} f(w) g(e^x) \frac{|\varepsilon|^{-1}}{w} \exp\left(-\frac{1+e^{2x}}{2w^{1/\varepsilon}}\right) I_{\sqrt{2z}}\left(\frac{e^x}{w^{1/\varepsilon}}\right) dx\, dw.$$

On the other hand, we have tautologically

$$L(z) = \int_0^\infty \int_{\mathbf{R}} f(w) g(e^x) \int_0^\infty e^{-zt} \frac{e^{-\frac{x^2}{2t}}}{\sqrt{2\pi t}} a_{\varepsilon,t}(x,w)\, dt\, dx\, dw$$

sucessively expressing the expectations defining this Laplace transform as obtained by integrating against the respective $\varepsilon$–power Asia densities and expressing these latter densities as conditional $\varepsilon$–power Asia densities $a_{\varepsilon,t}(x,w)$ integrated against the semigroups of the conditions $B_t = x$. Since $f$, $g$ are arbitrary maps, we get on comparison

$$\mathscr{L}\left(t \mapsto \frac{1}{\sqrt{2\pi t}} \exp\left(-\frac{x^2}{2t}\right) a_{\varepsilon,t}(x,w)\right)(z) = \frac{|\varepsilon|^{-1}}{w} \exp\left(-\frac{1+e^{2x}}{2w^{1/\varepsilon}}\right) I_{\sqrt{2z}}\left(\frac{e^x}{w^{1/\varepsilon}}\right).$$

for any real $z \geq z_0$. We extend the validity of this identity to the whole complex half–plane $H = \{z | \operatorname{Re}(z) \geq z_0\}$ using analytic continuation. For this we have to establish analyticity on $H$ of both sides of the identity considered as functions in $z$. Indeed, granting this for the moment, the identity theorem then applies and extends the validity of the above Laplace transform identity to the whole complex half–plane $H$ as desired

To establish the required analyticity properties, first consider the right hand side of the above identity. It is analytic on the whole right complex half–plane in particular since the Bessel function is entire in its degree and the square root is analytic on the complex plane with the non–positive reals deleted.

What regards analyticity in $z$ of the Laplace transform of the left hand side of the above identity, notice that an application of the *analyticity criterion* §7 Lemma reduces this to the integrability of its integrands $f(z,t) = \exp(-zt)(2\pi t)^{-\frac{1}{2}} \exp(-x^2/(2t)) a_{\varepsilon,t}(x,w)$ if $z$ is in any compact subset of $H$ and $t$ ranges over the positive real line. This, however, is implied by §5 Proposition, and the proof of the Lemma is complete.

**10. Analytic Laplace inversion for Yor type $\varepsilon$–power Asia densities:** Analytic Laplace inversion of the transforms for the conditional $\varepsilon$–power Asia densities $a_{\varepsilon,t}(x,w)$ of §9 Lemma is along the lines of of [**Y80**, §5]. In fact, it reduces to the explicit Laplace inversion of §6 Lemma, and the precise result is checked to be as follows

**Lemma:** *Let $\varepsilon$ be any non–zero real. For any reals $t > 0$ and $x$, the conditional $\varepsilon$–power Asia density at any $w > 0$ is given by the following contour integral*

$$\frac{1}{\sqrt{2\pi t}} \exp\left(-\frac{x^2}{2t}\right) a_{\varepsilon,t}(x,w) = \frac{b_{\varepsilon,t,w}}{2\pi i} \int_{\log C_{\theta,R}} e^{-\frac{\xi^2}{2t} + \frac{\exp(x)}{2w^{1/\varepsilon}} \cosh(\xi)} \sinh(\xi)\, d\xi,$$

*with real parameters $\theta$ in $[\frac{\pi}{2}, \pi]$ and $R \geq 1$.*

The factors $b_{\varepsilon,t,w}$ herein are given by

$$b_{\varepsilon,t,w} = \frac{|\varepsilon|^{-1}}{\sqrt{2\pi t}} e^x w^{-(1+\frac{1}{\varepsilon})} \exp\left(-\frac{1+e^{2x}}{2w^{1/\varepsilon}}\right)$$

and $\log_{\theta,R}$ are the contours of §2.

**11. Proof of the Yor type $\varepsilon$–power Asia densities:** For the proof of the Yor form of the $\varepsilon$–power Asia densities of §3 Theorem, we again compute in two ways the expectation

$$E(f) = E\left[f\left((A_t^{(\nu)})^\varepsilon\right)\right]$$

for any measurable function $f$ of the non–negative real line into itself. For the first of these, successively apply the Girsanov identity of §8 and formally express the resulting expectation by double integration against the appropriate densities. This gives

$$E(f) = e^{-\frac{1}{2}\nu^2 t} \int_{\mathbf{R}} \int_0^\infty e^{\nu x} f(w) \frac{1}{\sqrt{2\pi t}} e^{-\frac{x^2}{2t}} a_{\varepsilon,t}(x,w)\, dw\, dx\,.$$

On substitution for the density from §10 Lemma and the change of variables $y = \exp(x)$,

$$E(f) = \int_0^\infty \int_0^\infty y^\nu f(w) c^* e^{-\frac{y^2}{2w^{1/\varepsilon}}} \frac{1}{2\pi i} \int_D e^{-\frac{\xi^2}{2t} + \frac{y}{w^{1/\varepsilon}} \cosh(\xi)} \sinh(\xi)\, d\xi\, dw\, dy\,.$$

setting $c^* = c_{\nu,\varepsilon,t}(w)/\sqrt{(2w^{1/\varepsilon})^{\nu+1}}$ and $D = \log C_{\theta,R}$. Herein a Tonelli argument permits to interchange the order of the first two integrations. Changing variables $(2w^{\frac{1}{\varepsilon}})^{\frac{1}{2}} x = y$ in the resulting middle integral, we obtain

$$E(f) = \int_0^\infty f(w)\, c_{\nu,\varepsilon,t}(w) \int_0^\infty x^\nu e^{-x^2} \frac{1}{2\pi i} \int_D e^{-\frac{\xi^2}{2t}} \sinh(\xi)\, e^{\frac{2x \cosh(\xi)}{\sqrt{2w^{1/\varepsilon}}}}\, d\xi\, dx\, dw\,.$$

On the other hand, expressing $E(f)$ as obtained by integration against the $\varepsilon$–power Asia density of index $\nu$ at time $t$, we have tautologically

$$E(f) = \int_0^\infty f(w)\, \alpha_t^{(\nu,\varepsilon)}(w)\, dw.$$

The identity of §3 Theorem so follows on comparison, and the proof is complete.

# Part IV
# Probabilistic part of the proof of the
# Hermite function form of the $\varepsilon$–power Asia densities

**12.  Laplace transforms of the Hermite $\varepsilon$–power Asia densities:**  The proof of the Hermite function form of the $\varepsilon$–power Asia densities $\alpha^{(\nu,\varepsilon)}$ of §4 Theorem is by a Laplace transform analysis. As its first step, this section computes a second integral representation for the Laplace transform of these densities. The precise result is the following

**Lemma:**  *Let $\varepsilon$ be any non–zero real. If $\nu \geq 0$, we have the Laplace transform*

$$\mathscr{L}\left(t \mapsto \alpha_t^{(\nu,\varepsilon)}(w)\right)(z) = \frac{(2w^{\frac{1}{\varepsilon}})^{\frac{\nu}{2}}}{|\varepsilon|\, w} e^{\frac{-1}{2w^{1/\varepsilon}}} \int_0^\infty I_{\sqrt{2z+\nu^2}}\left(\frac{2x}{\sqrt{2w^{1/\varepsilon}}}\right) x^{\nu-1} e^{-x^2}\, dx$$

*for any positive reals $w$ and any complex $z$ with a sufficiently big positive real part.*

To prove this result let $f$ be any measurable function on the non–negative real line, and consider the map $F$ on the positive real line given by

$$F(t) = E\left[f\left((A_t^{(\nu)})^\varepsilon\right)\right].$$

The idea is to compute in two ways its Laplace transform

$$\mathscr{L}(F)(z) = \int_0^\infty e^{-zt} F(t)\, dt$$

with $z \geq z_0$ any sufficiently big real. First, §5 Proposition applied with $g = 1$ gives on changing variables $(2w^{\frac{1}{\varepsilon}})^{\frac{1}{2}} x = \rho$

$$\mathscr{L}(F)(z) = \int_0^\infty f(w)\, \frac{(2w^{\frac{1}{\varepsilon}})^{\frac{\nu}{2}}}{|\varepsilon|\, w} e^{\frac{-1}{2w^{1/\varepsilon}}} \int_0^\infty I_{\sqrt{2z+\nu^2}}\left(\frac{2x}{\sqrt{2w^{1/\varepsilon}}}\right) x^{\nu-1} e^{-x^2}\, dx\, dw$$

for $z \geq z_0$ any sufficiently big positive real. On the other hand, expressing in the defining integral of $\mathscr{L}(F)(z)$ the expectations as given by integration against the $\varepsilon$–power densities, Tonelli's theorem applies to give

$$\mathscr{L}(F)(z) = \int_0^\infty f(w)\, \mathscr{L}\left(\alpha_t^{(\nu,\varepsilon)}(w)\right)(z)\, dw.$$

Since $f$ is arbitrary, we have on comparison

$$\mathscr{L}\left(t \mapsto \alpha_t^{(\nu)}(w)\right)(z) = \frac{(2w^{\frac{1}{\varepsilon}})^{\frac{\nu}{2}}}{|\varepsilon| w} e^{\frac{-1}{2w^{1/\varepsilon}}} \int_0^\infty I_{\sqrt{2z+\nu^2}}\left(\frac{2x}{\sqrt{2w^{1/\varepsilon}}}\right) x^{\nu-1} e^{-x^2} dx$$

for any reals $z \geq z_0$. The validity of this identity extends to the whole complex half–plane $\{z | \text{Re}(z) \geq z_0\}$ by an analytic extension argument similar to that for §9 Lemma. The additional fact to be shown is that the analyticity properties of the right hand's side integrands in $z$ are preserved on integration. For this use the asymptotic behaviour of the Bessel functions near zero, see [**L**, §5.7], and towards infinity, see [**L**,§5.11], to construct as a first step functions majorizing the absolute values of these integrands which depend continuously on the parameter $z$. The desired analyticity is then is a consequence of the *analyticity criterion* §7 Lemma, and the proof of the Lemma is complete.

**13. Analytic Laplace inversion for the Hermite $\varepsilon$–power Asia densities:** As a second step in the Laplace transform analysis leading to the Hermite function form of Yor's Asia density $\alpha^{(\nu,\varepsilon)}$ of §4 Theorem, this section analytically inverts the Laplace transform of §12 Lemma. The precise result is the following

**Lemma:** *The assertions of §4 Theorem are valid if $\nu > 0$.*

The idea for the inversion of §12 Lemma's Laplace transform is to justify Laplace inversion under the integral sign and apply the inversion result of §6 Lemma. For this again express the Bessel function factor as the Hankel contour integral from [**WW**, 17·231, p.362]:

$$I_\mu(\eta) = \frac{1}{2\pi i} \int_D e^{-\mu \cdot \xi + \eta \cdot \cosh(\xi)} d\xi,$$

with $\mu = (2z+\nu^2)^{1/2}$, with $\eta = 2x(2w^{\frac{1}{\varepsilon}})^{-\frac{1}{2}}$, and with the contour $D = \log C_{\theta,R}$. Then study the asymptotic behaviour of the resulting double integral's integrand for $(z,x,\xi)$ to infinity. Hereby $z$ varies on a parallel to the imaginary axis suffciently deep within the right half plane, $x$ varis in the positive real line, and $\xi$ varies in the contour $D$. First restrict to contours with $\theta$ in $(\frac{\pi}{2}, \pi]$ and $R$ so large that $D$ and its image under squaring are contained in the right half plane. As in the proof of §6 Lemma, the above integrand then is seen to be exponentially decreasing to zero with such $(z,x,\xi)$ to infinity. Laplace inversion under the integral sign is thus justified in particular. Using §6 Lemma,

$$\alpha_t^{(\nu,\varepsilon)}(w) = c_{\nu,\varepsilon,t}(w) \int_0^\infty e^{-x^2} x^\nu \frac{1}{2\pi i} \int_D e^{-\frac{\xi^2}{2t} + \frac{2x\cosh(\xi)}{\sqrt{2w^{1/\varepsilon}}}} \sinh(\xi) d\xi dx$$

on collecting constants. The asymptotic behaviour established above moreover implies that the order of the two integrations can be interchanged. Since $\nu$ is positive, the resulting inner integral is a value of a certain Hermite function of negative degree recalling the latter's integral representation for negative degrees from §2. As precise result we get

$$\alpha_t^{(\nu,\varepsilon)}(w) = \frac{c_{\nu,\varepsilon,t}(w)}{2\pi i} \int_{\log C_{\theta,R}} e^{-\frac{\xi^2}{2t}} \sinh(\xi)\, \Gamma(\nu+1) H_{-(\nu+1)}\!\left(-\frac{\cosh(\xi)}{\sqrt{2w^{1/\varepsilon}}}\right) d\xi.$$

To lift the restrictions on $\theta$ and $R$, notice that the integrand of this last integral is entire as a function of $\xi$. Using the asymptotic expansion for Hermite functions

$$H_\mu(z) = (2z)^\mu \sum_{k=0}^{n-1} \frac{(-\mu)_{2k}}{k!} \cdot \frac{(-1)^k}{(2z)^{2k}} + O\left(\frac{1}{|z|^{2n-\operatorname{Re}(\mu)}}\right)$$

of [**L**, 10.6.7, p.292], valid for any complex $z$ with $|\arg(z)| < 3\pi/4$, the integrals are finite if $\theta$ is in $[\frac{\pi}{2}, \pi]$ and $R \geq 1$. Cauchy's theorem applies and shows that the integral is independent of $\theta$ and $R > 1$. This extends to constellations with $R = 1$ by continuity, and the proof of the Lemma is complete.

## Part V      Analytical preliminaries

**14. Analyticity properties of Yor type $\varepsilon$–power Asia densities:** This section considers the contour integral factors of the Yor type form of the $\varepsilon$–power densities of §3 Theorem as functions in the index. More precisely, for any non–zero real $\varepsilon$, any positive reals $t$ and $w$, we thus consider the functions which for any complex $\mu$ are given by

$$P_{\varepsilon,t,w}(\mu) = \int_0^\infty x^\mu e^{-x^2} \psi_{\varepsilon,t,w}(x) \, dx,$$

where $\psi_{\varepsilon,t,w}$ are the functions on the real line given by the integrals

$$\psi_{\varepsilon,t,w}(x) = \frac{1}{2\pi i} \int_{\log C_{\theta,R}} e^{-\frac{\xi^2}{2t}} \sinh(\xi) \, e^{\frac{2x\cosh(\xi)}{\sqrt{2w^{1/\varepsilon}}}} \, d\xi$$

over any logarithmicalized Hankel contour $\log C_{\theta,R}$ of §2. The result to be proved is the following

**Proposition:** *All functions $P_{\varepsilon,t,w}$ are entire.*

**Remark:** For any real $\nu$, we have $\alpha_t^{(\nu,\varepsilon)}(w) = c_{\nu,\varepsilon,t}(w) P_{\varepsilon,t,w}(\nu)$.

With the Remark being clear by construction, notice that the proof of the Proposition reduces to a local question. In fact, dropping references to indices $\varepsilon$, $t$, $w$ and $k$ to simplify the notation, it is sufficient to prove entireness of $P$ for $\mu$ ranging within any compact subset $V$ of the complex plane with non–empty interior. As an application of the *analyticity criterion* §7 Lemma this further reduces to construct an integrable function $g$ on the positive real line such that

$$|p(\mu, x)| \leq g(x)$$

for any $x > 0$ and $\mu$ in $V$ setting $p(\mu, x) = \exp(-x^2) x^\mu \psi(x)$. For this construction first apply Cauchy's theorem to see that $\psi$ is independent of the contours $\log C_{\theta,R}$ chosen. The function $\psi$ that results from taking $\theta = \pi$ and $R = 1$ is checked to coincide with the one considered in [**Y**, 6.g], whence

$$\lim_{x \downarrow 0} x^\mu \psi(x) = 0$$

in particular for any $\mu$ in $V$. For any $c > 0$, sending $\mu$ in $V$ and $x$ in $(0, c]$ to $x^\mu \psi(x)$ thus extends to a continuous map on $V \times [0, c]$. Let $C_V$ be the maximum its absolute value takes. On the interval $[c, \infty)$ majorize $|\psi|$ by a positive constant $D_c$ depending on $c$. For instance, choosing $\theta = \pi$ and $R = 1$, take for this

$$D_c = \frac{e^{\frac{\pi^2}{2t}}}{\pi} \int_0^\infty e^{-\frac{y^2}{2t}} \sinh(y) \, e^{-\frac{2c \cosh(y)}{\sqrt{2w^{1/\varepsilon}}}} \, dy \, .$$

To continue, let $\pi_{\mu_V}$ denote the map on the positive real line sending any $x$ to its $\mu_V$–th power where $\mu_V$ is the maximum of the real parts of all elements in $V$. Moreover choosing $c > 1$, the function $g$ on the positive real line defined by

$$g(x) = \exp(-x^2)(C_V \mathbf{1}_{[0,c)} + D_c \pi_{\mu_V} \mathbf{1}_{[c,\infty)})(x)$$

then has the desired majorizing properties, and the proof of the Proposition is complete.

**15.  A connection with Hermite functions:**  As a first explanation for the relatedness of the two forms of the $\varepsilon$–power Asia densities this section interprets the integrals $P_{\varepsilon,t,w}(\mu)$ of §14 in terms of the Hermite functions $H_\mu$ recalled in §2. More precisely, consider for any non–zero real $\varepsilon$ and any positive reals $t$ and $w$ the functions

$$F_{\varepsilon,t,w}(\mu) = \frac{1}{2\pi i} \int_{\log C_{\theta,R}} f_{\varepsilon,t,w}(\mu, \xi) \, d\xi$$

setting

$$f_{\varepsilon,t,w}(\mu, \xi) = e^{-\frac{\xi^2}{2t}} \sinh(\xi) \, H_{-(\mu+1)}\!\left(-\frac{\cosh(\xi)}{\sqrt{2w^{1/\varepsilon}}}\right) d\xi \, .$$

Then the precise result is the following

**Lemma:**  *For any complex number $\mu$ with $\mathrm{Re}\,(\mu) > -1$, we have*

$$P_{\varepsilon,t,w}(\mu) = \Gamma(\mu+1) \, F_{\varepsilon,t,w}(\mu)$$

*for any positive reals $t$ and $w$ and non–negative integers $k$.*

To see this revert to definitions and write out the integrals $P_{\varepsilon,t,w}(\mu)$ of §14 as double integrals. The Lemma's identity follows on interchanging the order of integration as soon as we know the integrability of this double integral's integrand. This however is shown by an argument analogous to that for §6 Lemma or §13 Lemma, and the proof is complete.

**16. Analyticity results:** This section adresses analyticity of the contour integrals $F_{\varepsilon,t,w,k}$ of §15 as functions in the complex variable $\mu$. Thus recall for any non–zero real $\varepsilon$ and any positive reals $t$ and $w$ the functions

$$F_{\varepsilon,t,w}(\mu) = \frac{1}{2\pi i} \int_{\log C_{\theta,R}} f_{\varepsilon,t,w}(\mu,\xi)\, d\xi$$

setting

$$f_{\varepsilon,t,w}(\mu,\xi) = e^{-\frac{\xi^2}{2t}} \sinh(\xi)\, H_{-(\mu+1)}\!\left(-\frac{\cosh(\xi)}{\sqrt{2}w^{1/\varepsilon}}\right) d\xi$$

which via §14 Remark and §15 Lemma also determine the analyticity properties of the $\varepsilon$–power Asia densities. The result to be proved is the following

**Proposition:** *All functions $F_{\varepsilon,t,w}$ are entire.*

The proof of the Proposition is similar to that of §14 Proposition. First, establishing this result is a local question. In fact, dropping reference to indices $\varepsilon$, $t$ and $w$ to simplify the notation, it is sufficient to prove entireness of $F$ for $\mu$ ranging within any compact subset $V$ of the complex plane with non–empty interior. Since in particular Hermite functions are entire functions in their degree, an application of the *analyticity criterion* §7 Lemma, again reduces this to show that there is an integrable function $g$ on $\log C_{\theta,R}$ with $|f(\mu,\xi)| \le g(\xi)$ for any $\xi$ in $\log C_{\theta,R}$ and $\mu$ in $V$. The construction of such a function is based on certain properties of Hermite functions. Using the recurrence rule for Hermite functions

$$H_{\mu+1}(z) - 2z \cdot H_\mu(z) + 2\mu \cdot H_{\mu-1}(z) = 0$$

of [**L**, 10.4.6, p.289] first reduce to $V$ being contained in the left complex half plane. Then use the uniform asymptotic expansion of Hermite functions

$$H_\mu(z) = (2z)^\mu \sum_{k=0}^{n-1} \frac{(-\mu)_{2k}}{k!} \cdot \frac{(-1)^k}{(2z)^{2k}} + \mathrm{O}\!\left(\frac{1}{|z|^{2n-\mathrm{Re}\,(\mu)}}\right)$$

of [**L**, 10.6.7, p.292], valid for any complex $z$ with $|\arg(z)| < 3\pi/4$, to see that the factors $\exp(-\xi^2/2t)$ dominate the asymptotic behaviour of the integrands of $F$ with $|\xi|$ to infinity. Moreover using the compactness of $V$ thus construct a majorizing function $g$ for the integrands depending on $V$ but independent of the single elements $\mu$ in $V$. The proof of the Proposition is complete.

# Part VI
# Analytic part of the proof of the
# Hermite function form of the Asia density

**17. Coincidence of the integrals:** This section relates the two integrals for the $\varepsilon$–power Asia density of §3 Theorem and §4 Theorem for processes of general indices $\nu$. More precisely, for any non–zero reals $\varepsilon$, and any positive reals $t$ and $w$, we have on the one hand the functions $P_{\varepsilon,t,w}$ for the Yor type integral discussed in §14, and on the other hand the functions $F_{\varepsilon,t,w}$ of §15 for the Hermite function form integral. The main result is then the following generalization of §15 Lemma.

**Proposition:** *For any positive reals $t$ and $w$, the functions $\mu \mapsto \Gamma(\mu+1)\, F_{\varepsilon,t,w}(\mu)$ and $\mu \mapsto P_{\varepsilon,t,w}(\mu)$ are entire, and we have*

$$P_{\varepsilon,t,w}(\mu) = \Gamma(\mu+1)\, F_{\varepsilon,t,w}(\mu)$$

*for any complex $\mu$.*

Indeed, for any fixed positive reals $t$ and $w$, the function $\mu \mapsto F_{\varepsilon,t,w}(\mu)$ is entire from §16 Proposition. The function $\mu \mapsto P_{\varepsilon,t,w}(\mu)$ has been proved entire in §14 Proposition. From §15 Lemma we have the equality

$$P_{\varepsilon,t,w}(\mu) = \Gamma(\mu+1)\, F_{\varepsilon,t,w}(\mu)$$

in particular for any positive $\mu$. Now apply the identity theorem. This equality then holds for any complex $\mu$ first in the sense of meromorphic functions since its gamma function factor has simple poles in the negative integers. With its left hand side an entire function, the right hand side of this identity cannot have any poles, and so this equality holds for any complex number $\mu$ between entire functions. This completes the proof of the Proposition.

**18. Proof of the Hermite function form of the Asia density:** To give a proof of the Hermite function form of the Asia density of §4 Theorem, recall there are two forms of the $\varepsilon$–power Asia densities $\alpha^{(\nu,\varepsilon)}$. Using the concepts of §2, we have on the one hand from §3 Theorem for any real $\nu$ the Yor type integral representation

$$\alpha_t^{(\nu,\varepsilon)}(w) = c_{\varepsilon,\nu,t}(w)\, P_{\varepsilon,t,w}(\nu)$$

for any positive reals $t$ and $w$, and with $P_{\varepsilon,t,w}$ the function discussed in §14. On the other hand, for any positive $\nu$, we have from §13 Lemma the Hermite function form of the density

$$\alpha_t^{(\nu,\varepsilon)}(w) = c_{\varepsilon,\nu,t}(w)\, \Gamma(\nu+1)\, F_{\varepsilon,t,w}(\nu)$$

for any positive reals $t$ and $w$. It so remains to extend this last identity to all real $\nu$.

For this first notice that from their definition in §2 the factors $c_{\varepsilon,\nu,t}(w)$ are entire functions in $\nu$ in particular. Using the above Yor type representation of the power densities, for any fixed positive reals $t$ and $w$, the function $\nu \mapsto c_{\varepsilon,\nu,t}(w) P_{\varepsilon,t,w}(\nu)$ represents the Asia density for all real $\nu$. Using §17 Proposition, we have the identity between entire functions

$$P_{\varepsilon,t,w}(\nu) = \Gamma(\nu+1)\, F_{\varepsilon,t,w}(\nu)$$

in particular for any real $\nu$. Thus $\nu \mapsto c_{\varepsilon,\nu,t}(w)\Gamma(\nu+1)\, F_{\varepsilon,t,w}(\nu)$ represents the $\varepsilon$–power Asia density for all real $\nu$. And the precise meaning of this last statement is that the above identity for the Hermite function form of the Asia density is valid for all real $\nu$, as it was to be shown. The proof of §4 Theorem is complete.

**19. Proof of §4 Remark:** The proof of §4 Remark is based on the following observation. Since $\nu \mapsto c_{\varepsilon,\nu,t}(w)\Gamma(\nu+1)\,F_{\varepsilon,t,w}(\nu)$ is entire from §17 Proposition, the poles of the gamma function factor are cancelled by zeros of the remaining factors. On inspection only the integral factor $\nu \mapsto F_{\varepsilon,t,w}(\nu)$ can develop such zeros, and the Remark follows.

**Acknowledgements:** Support by the *Deutsche Forschungsgemeinschaft* and hospitality of the *Institut für Mathematik* of the *Universität Mannheim* are gratefully acknowledged.

Author's address: Keplerstrasse 30, D-69469 Weinheim (Bergstrasse), Germany